\documentclass[10pt]{article}
\usepackage{amssymb}
\usepackage{amsthm}
\usepackage{graphicx}
\usepackage{graphics}
\usepackage{amsthm}

\input pictex.tex

\newcommand{\grad}{\mathrm{o}}

\theoremstyle{definition}

\input epsf

\begin{document}

\date{}
\author{Andr\'es Navas}

\title{The amazing story of a forgotten golden flag}
\maketitle

\begin{small}
\noindent{\bf Abstract.} We describe the most probable geometric design of the Chilean Independence Flag, 
which uses the golden ratio in many of its components. We also discuss some related historical aspects.
\end{small}

\vspace{0.4cm}

Several national flags incorporate mathematical elements in their design, which in most cases are fixed by law and 
sometimes are even explicitly stated in the corresponding national constitution. This is for instance the case of the flag 
of the United States of America, for which a nice and careful study about the symmetry rules for the arrangement of the stars 
is given in \cite{KT}. Most frequently, the involved elements are of a geometric nature, and are  
related to the proportions of the different components. Certainly, the most remarkable example is that of Nepal's flag, 
the only non rectangular national flag \cite{W1}, for which a quite intricate geometric construction \cite{Ba} 
yields the width-height ratio 

$$\frac{24 + \frac{297 - 180 \sqrt{2}}{92 - 36 \sqrt{2}} \left( 1 + \frac{8 - 3\sqrt{2}}{\sqrt{118 - 48\sqrt{2}} -6} \right) }
{32 + \frac{297 - 180 \sqrt{2}}{92 - 36 \sqrt{2}} \left( 1 + \frac{6}{(8 - 3\sqrt{2}) \big( \sqrt{1 + \frac{18}{41 - 24\sqrt{2}}} - 1\big)} \right)} 
\sim 0.820...$$
In many cases, these proportions are in relation to the golden mean, either to its Fibonacci approximations (e.g. Palau, Poland and Sweden's 
flags \cite{W1}) or its genuine value 
(e.g. Togo's flag \cite{W3}). The purpose of this Note is not at all to give a complete account of all mathematical aspects of the 
different national flags. Instead, we would like to draw the attention to the geometric beauty of another one, the Chilean Independence 
Flag, whose design  
is also based on golden proportions and remains widely unknown despite its almost two hundred years of history.

\begin{figure}[ht!]
\centering
\includegraphics[width=26.5mm]{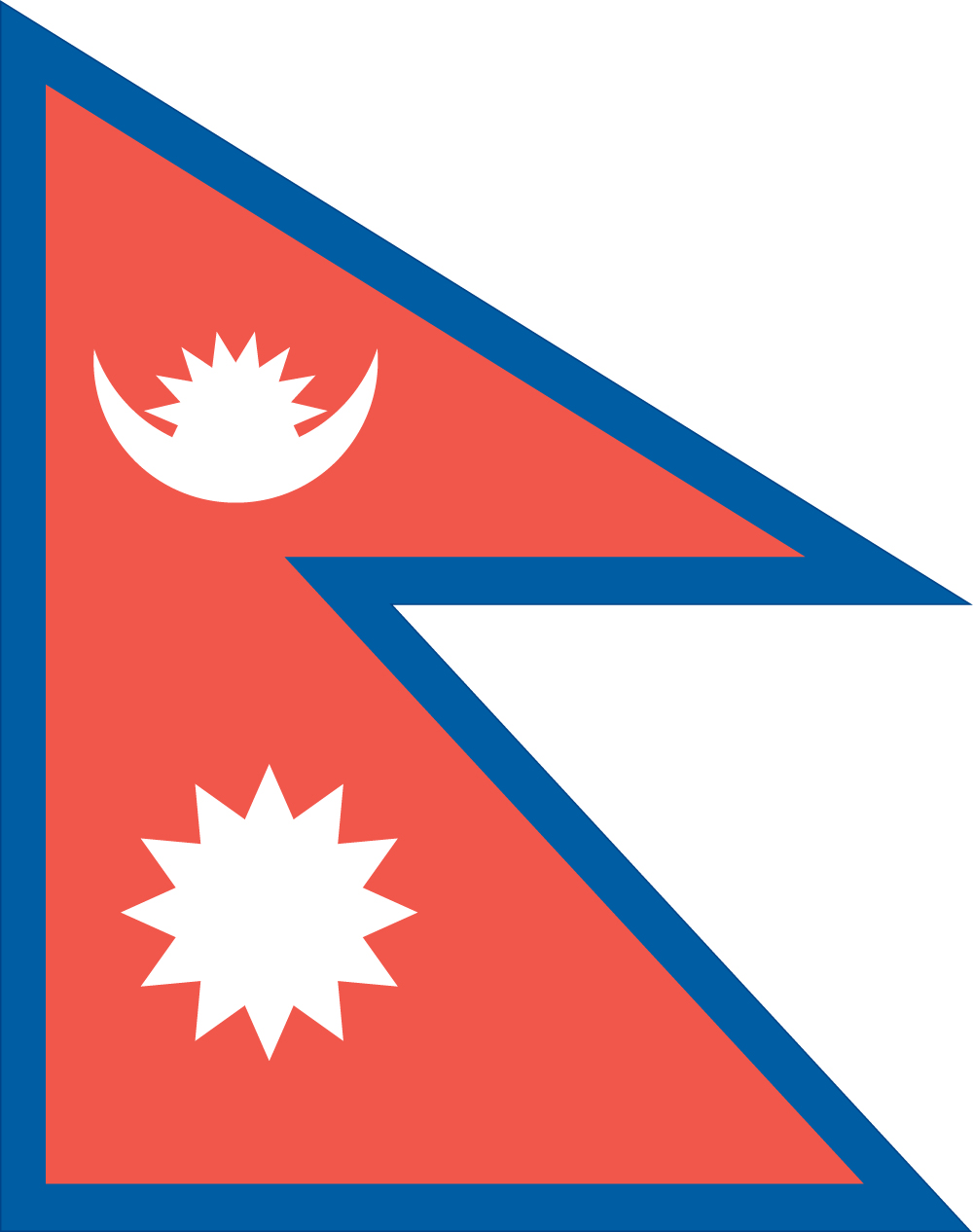}
\qquad\qquad
\includegraphics[width=55mm]{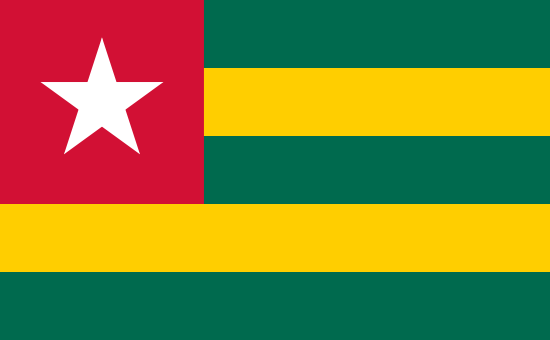}
\caption{On the left, Nepal's flag; on the right, Togo's flag (for which the width-height ratio equals the golden mean). \label{nepal-togo}}
\end{figure}

As it was the case of most independence movements in Latin America in the early $19^{th}$ century, the Chilean one was 
guided by the new political and social ideas that inspired prior processes in both the United States of America and France. 
Besides, the Chilean process was led by a masonic intelligentsia (the so-called {\em Logia Lautaro}\footnote{{\em Lautaro}  
was, in the $16^{th}$ century, the most important military leader (the {\em Cacique}) of the {\em Mapuche}, the largest 
indian community in the country that resisted the Spanish colonization.}) that transferred its symbolism to many of the 
national emblems \cite{B}. For instance,  the colors used in the current Chilean flag are white, blue and red. 

Although the original Independence Flag was quite similar to it, there are several remarkable differences. For example, 
following the explicit order of the military leader Bernardo O'Higgins, the Independence Flag included two shields (one with 
an obelisk and another one with a volcano), as well as an octogonal star (inside a pentagonal one) representing the planet Venus, 
as was usually depicted by the Mapuche. However, from a mathematical point of view,  the most 
spectacular differences concern geometry, as explained below.

It is worth emphasizing that no official document concerning the design of the Independence Flag remains. 
Actually, it is not even clear that such a document ever existed. In particular, the author of this design is unknown, although 
most historians claim that it is Antonio Arcos, a Spanish military engineer who fought on the Chilean side during the Independence War. 
Besides, there is just one copy of the flag, which was indeed used during the Independence Signature Ceremony (February $12^{th}$, 1818) 
and is preserved in the National Historic Museum at Santiago.\footnote{This 
has been so since the independence, except that in 1980 it was stolen as a sign of resistance to the Pinochet 
dictatorship. (It was belatedly returned in 2003.)} Since then, only one person -- the musician and philosopher 
Gast\'on Soublette \cite{S} -- noticed that the design 
incorporates the use of the golden mean. Unfortunately, his analysis was incorrect as it was guided mainly by symbolic elements 
and not by geometric ones.

What follows is a summary of the analysis in \cite{Na} leading to 
what most probably was the geometrical design of the Chilean Independence Flag. Notice that the description 
below is modeled on, yet not always coincides with, the measures of the existing flag, as this is somewhat 
deformed and has been restored in many ways, without paying much attention to its proportions; see \cite{MRCM}.  

\vspace{-0.1cm}

\begin{figure}[ht!]
\centering
\includegraphics[width=102mm]{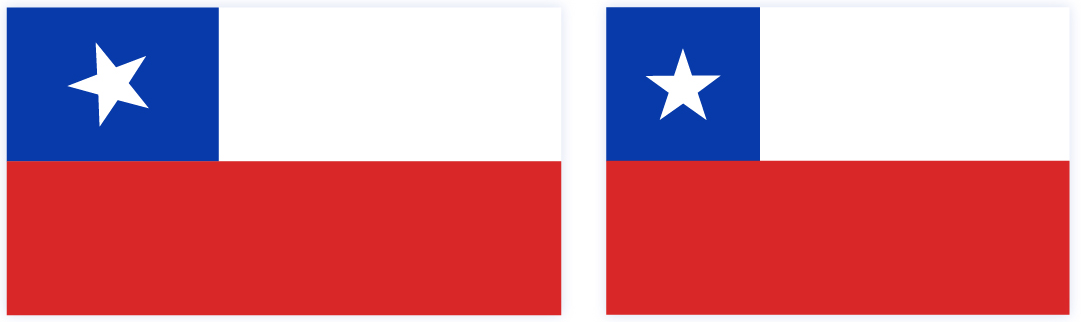}
\vspace{-0.2cm}
\caption{On the left, the (reconstructed) Independence Flag (without the shields and the octogonal star); 
on the right, the current flag. \label{flags}}
\end{figure}

Like the current flag, the Independence Flag was composed of three regions: a red rectangle at the bottom, and a blue and a white 
rectangle at the top (the white one being to the right of the blue one), all of them with the same height. The ratio of the widths of the white and 
blue parts equals the golden mean. The blue part corresponds to a rectangle for which the proportion between the 
height and width equals 
$$\mathrm{tan} (36^{\grad}) 
= \frac{\sqrt{10 - 2 \sqrt{5}}}{1 + \sqrt{5}} 
= \frac{\sqrt[4]{5}}{\sqrt{2+\sqrt{5}}} \sim 0.726...$$ 
This choice leads to the quite elegant angle configuration below, in which several 
isosceles triangles of angles $36^{\grad}, 72^{\grad},72^{\grad}$ appear, so that golden proportions are omnipresent.\footnote{This is 
to be compared with the Chinese national flag, whose configuration also incorporates angles of $36^{\grad}$; see \cite{china}.} Besides, this 
makes the width-height ratio of the flag to be equal to 
$$\frac{2 + \sqrt{5}}{\sqrt{10 - 2 \sqrt{5}}} \sim 1.801...$$
Finally, a white pentagonal star is depicted in the blue rectangle with its center in the intersection point of the diagonals, as shown below. 
The size of this star is very special: the ratio between the height of the rectangle and the diameter of 
the star's circumcircle equals, again, the golden mean.

\begin{figure}[ht!]
\centering
\includegraphics[width=61mm]{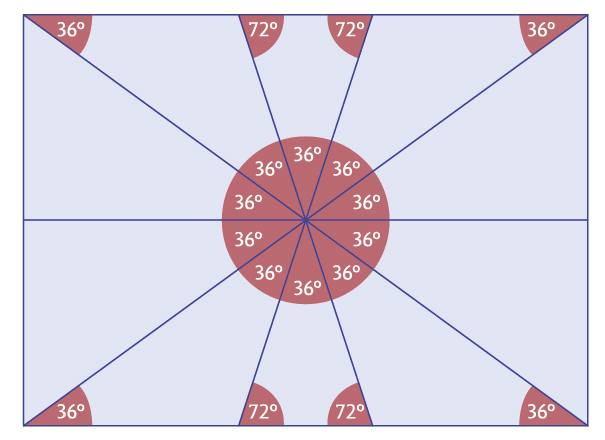}
\qquad
\includegraphics[width=61mm]{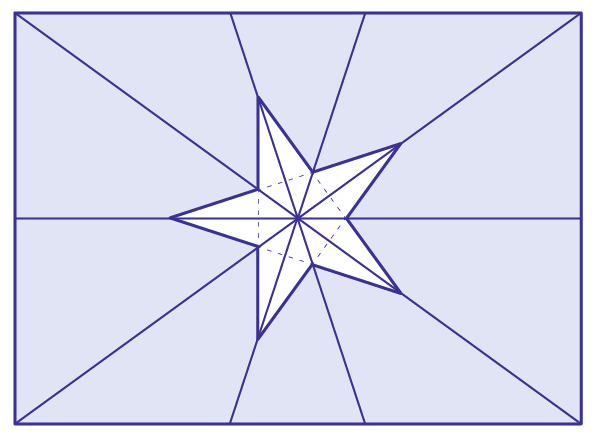}
\vspace{-0.1cm}
\caption{The angle configuration and the star. \label{angles}}
\end{figure}

\newpage

With this configuration, the Independence Flag was elaborated with a 2.4 m width size. Concerning the current flag, it was already in use 
in the middle of the $19^{th}$ century (yet the final law establishing its proportions dates back to 1912): the width-height proportion is $3 \! : \! 2$, 
so that it is composed of six squares (3 red squares, 2 white ones and a blue one); the star is centered at the center of the blue square, and 
the diameter of its circumcircle is half the length of the side of the square. Needless to say, this was certainly a regrettable simplification. 

It is now the task of the historians to determine 
the reasons for this very quick replacement and the lack of major information about the original Independence Flag:\ was this due 
to the complexity of the design and the difficulties for the elaboration at that time, 
or was it the result of internal disputes between the independence leaders and the first rulers of the 
country?\footnote{Two of the three independence leaders, Jos\'e Miguel Carrera and Manuel Rodr\'iguez, were executed, most probably with the 
consent of the Logia Lautaro; the third one, O'Higgins, ruled the country up to 1823, when he was exiled to Per\'u, where he died in 1842.}

\begin{figure}[ht!]
\centering
\includegraphics[width=85mm]{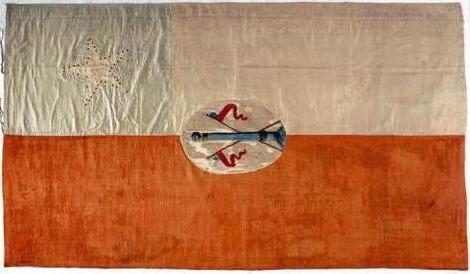}
\vspace{-0.18cm}
\caption{A picture of the genuine Independence Flag. \label{flags}}
\end{figure}


\vspace{0.01cm}

\noindent{\bf Acknowledgments.} I would like to thank Fanny Espinoza, as well as all people in the National Historic Museum at Santiago, 
for the reception and the permission to examine the Independence Flag. Also, I would like to thank Hern\'an Aburto for introducing me to 
the subject, Carolina Mu$\tilde\mathrm{n}$oz for drawing some of the pictures of this Note, and Mahsa Allahbakhshi, Eduardo Friedman, Rodrigo Hern\'andez,  Godofredo Iommi, August J. Krueger, Mario Ponce, Zoran \v Suni\'c and Jessica Vidal for their many suggestions 
and corrections to the original version.

This work was funded by the Center of Dynamical Systems and Related Fields (Anillo Research Project 1103, CONICYT).


\begin{small}

\vspace{0.1cm}

\noindent Andr\'es Navas\\ 

\noindent Dpto. de Matem\'atica y Ciencia de la Computaci\'on, Universidad de Santiago de Chile\\ 

\noindent Alameda 3363, Santiago, Chile\\ 

\noindent email: andres.navas@usach.cl

\end{small}

\end{document}